\documentclass{aptpub}
\authornames{Jean-Marc Bardet, Paul Doukhan, Gabriel Lang, and Nicolas Ragache}
\shorttitle{Dependent Lindeberg central limit theorem}
\usepackage[english]{babel}
\usepackage{latexsym}
\usepackage{amsmath}
\usepackage{amssymb}
\usepackage{graphicx,epsfig}
\def\1{\mbox{1\hspace{-.35em}1}} 
\def\R{\mathbb{R}}

\def\N{\mathbb{N}}
\def\P{\mathbb{P}}
\def\E{\mathbb{E}}

\def\R{\mathbb{R}}

\def\Z{\mathbb{Z}}

\def\v{\mbox{Var\,}}
\def\lip{\mbox{Lip\,}}
\def\Lip{\mbox{Lip\,}}
\def\cov{\mbox{Cov}}%

\def\limitek{\renewcommand{\arraystretch}{0.5}
\begin{array}[t]{c}\stackrel{}{\longrightarrow} \\
{\scriptstyle k\rightarrow \infty} \end{array}
\renewcommand{\arraystretch}{1}}

\def\limiten{\renewcommand{\arraystretch}{0.5}
\begin{array}[t]{c}\stackrel{}{\longrightarrow} \\
{\scriptstyle n\rightarrow \infty} \end{array}
\renewcommand{\arraystretch}{1}}

\def\limiter{\renewcommand{\arraystretch}{0.5}
\begin{array}[t]{c}\stackrel{}{\longrightarrow} \\
{\scriptstyle r\rightarrow \infty} \end{array}
\renewcommand{\arraystretch}{1}}

\def\limiteloik{\renewcommand{\arraystretch}{0.5}
\begin{array}[t]{c}\stackrel{{\cal D}}{\longrightarrow} \\
{\scriptstyle k\rightarrow \infty} \end{array}
\renewcommand{\arraystretch}{1}}

\def\limiteloin{\renewcommand{\arraystretch}{0.5}
\begin{array}[t]{c}\stackrel{{\cal D}}{\longrightarrow} \\
{\scriptstyle n\rightarrow \infty} \end{array}
\renewcommand{\arraystretch}{1}}

\begin{document}
\title{Dependent Lindeberg central limit theorem \\and
some applications}

\authorone[Samos-Matisse-CES, Universit\'e Paris 1, CNRS UMR 8174]{Jean-Marc Bardet}
\authortwo[CREST and Samos-Matisse-CES, Universit\'e Paris 1, CNRS UMR 8174]{Paul Doukhan}
\authorthree[ENGREF, UMR MIA 518, INRA ENGREF INAP-G, ParisTech]{Gabriel Lang}
\authortwo[CREST]{Nicolas Ragache}

\addressone{Samos-Matisse-CES, Universit\'e Panth\'eon-Sorbonne, 90 rue de Tolbiac, 75013 Paris, FRANCE}
\addresstwo{LS-CREST, Timbre J340, 3 avenue Pierre Larousse, 92240 Malakoff, FRANCE}
\addressthree{ENGREF, 19, avenue du Maine, 75732 Paris Cedex 15 FRANCE}
\begin{abstract}
In this paper, a very useful lemma (in two versions) is proved: it
simplifies notably the essential step to  establish a Lindeberg
central limit theorem for dependent processes. Then, applying this
lemma to weakly dependent processes introduced in Doukhan and
Louhichi (1999), a new central limit theorem is obtained for
sample mean or kernel density estimator. Moreover, by using the
subsampling, extensions under weaker assumptions of these central
limit theorems are provided. All the usual causal or non causal
time series: Gaussian, associated, linear, ARCH($\infty$),
bilinear, Volterra processes,$\ldots$, enter this frame.\\
\end{abstract}



\section{Introduction}
This paper adresses the problem of the central limit theorem
(C.L.T.) for weakly dependent sequences with the point of view of
the classical Lindeberg method (see Petrov, 1995, for references).
For establishing a C.L.T. for a sequence $(S_n)_{n\in \N^*}$ of
random vectors (r.v.), a convenient and efficient method
(so-called the "Lindeberg's method" in the sequel) consists on
proving that for all functions $f$ with bounded and continuous
partial derivatives up to order $3$,
\begin{eqnarray}\label{method_lindeberg}
\Big|\E\big(f(S_n)-f(N)\big)\Big| \limiten 0,
\end{eqnarray}
with $N$ a Gaussian random variable not depending on $n$. Assume
that $S_n=X_1+\cdots+X_n$ where $(X_k)_{k\in \N^*}$ is a zero-mean
sequence with $(2+\delta)$-order finite moments for some
$\delta>0$, and consider $(Y_k)_{k\in \N^*}$ a sequence of
independent zero mean Gaussian r.v. such that the variance of
$Y_k$ and $X_k$ are the same. In order to obtain
(\ref{method_lindeberg}), we first show that this convergence is
satisfied when the sum of two terms converges to zero. The first
term is the sum of the $(2+\delta)$-order moments of $(X_k)_{1\leq
k \leq n}$. The second term is a sum of covariances between
functions of $(X_k)_{1\leq k \leq n}$ and $(Y_k)_{1\leq k \leq n}$
and reflects all the dependence structure. Three cases are thus
detailed in three different lemmata: first, the case of
independent r.v.s (here, the second term vanishes), then the
dependent case with general functions $f$, and finally the
dependent case with characteristic
functions that yields a very simple expression. \\
~\\
For applications of those lemmata, the class of weakly dependent
processes, introduced by Doukhan and Louhichi (1999), is selected
here. Roughly speaking, a process $X=(X_k)_{k\in \N^*}$ is said to
be a weakly dependent process if the covariance of any bounded and
Lipschitz function of "past" data of $X$ by any bounded and
Lipschitz function of "future" data of $X$ tends to $0$ when the
lag between the future and the past increases to $\infty$ (see a
more precise definition below).

\begin{center} Why should we use such dependence structures (instead {\it
e.g.} mixing)? \end{center}

Two main reasons motivate this choice.  Firstly, weak dependence
is a very general property including certain non-mixing processes:
{\it e.g.} Andrews (1984) explicited the simple example of an
autogressive process with Bernoulli innovations and proved that
such a model is not mixing in the sense of Rosenblatt (see for
instance Doukhan, 1994, or Rio, 2000, for references) while
Doukhan and Louhichi (1999) proved that such a process is weakly
dependent. More generally, under weak conditions, all the usual
causal or non causal time series are weakly dependent processes:
this is the case for instance of Gaussian, associated, linear,
ARCH($\infty$), bilinear, Volterra, infinite memory
processes,$\ldots$. Secondly, the dependence property is obtained
from the convergence to zero of covariances of the process (see
above). The second term in our lemmata writes as a sum of
covariances; weak dependence is therefore particularly accurate to
bound this term  (which is the essential step for proving the
Lindeberg C.L.T.). \\
~\\
Different applications of the lemmata are then presented for
weakly dependent processes. First, a C.L.T. for sample means is
established in Doukhan and Wintenberger (2006); in addition to the
previous lemma for characteristic function, the Bernstein block
method is required (in such a case, Bulinski and Shashkin, 2004
and 2005 used also this method; an alternative method is derived
in Rio, 2000, Coulon-Prieur and Doukhan, 2000, or Neumann and
Paparoditis, 2005). For weakly dependent processes, a C.L.T. for
the subsample mean is derived here directly from our lemmata. By
this way, the conditions required for such a theorem are weaker
than those required for the C.L.T. for the sample mean. For
instance, the subsampling of a long range dependent process
provides a C.L.T. for its sample mean, which is interesting for
obtaining confidence intervals or semi-parametric
tests (even if a large part of the sample is not used). \\
Finally, an application of the Lindeberg method to the kernel
density estimation is also given. By this way, the C.L.T. is
established under the same conditions (but with a more simple a nd
genral method) than in Coulon-Prieur and Doukhan (2000) for causal
processes. Its extension to non-causal processes is also proposed
here. The required conditions are of a different nature that the
usual conditions under strong mixing (see Robinson, 1983); but on
some examples of time series (for instance, the causal linear
processes) they imply the same decay rates of the coefficients. On
other examples (some non-causal time series), it is very difficult
to check the strong mixing property. Therefore the C.L.T. we
proved concerns a lot of new models. Moreover, a version of this
C.L.T. for subsampled kernel density estimator is given. Once
again, this allows to obtain a C.L.T. under weaker conditions.
With an adapted subsampling step, the asymptotic normality of this
estimator is established even in the case of long memory
processes, which provides usual confidence intervals on the
density or goodness-of-fit tests. \\
~\\
The paper is organized as follows. In Section \ref{lindsect}, the
Lindeberg method is presented. The Section \ref{sec3} is devoted
to a presentation of weakly dependent processes. Section
\ref{appli} contains different applications of the Lindeberg
method for weakly dependent processes while the proofs of the
different results are in the Section \ref{proofs}.
\section{Lindeberg method}\label{lindsect}
Let $(X_i)_{i\in \N}$ be a sequence of zero mean r.v. with values
in $\R^d$ (equipped with the Euclidean norm
$\|X_i\|^2=X_{i,1}^2+\cdots+X_{i,d}^2$ for
$X_i=(X_{i,1},\ldots,X_{i,d})$). Moreover, all along this paper we
will assume that $(X_i)_{i\in \N}$ satisfies,\\
~\\
{\bf Assumption $H_\delta$:} It exists $0<\delta\le1$ such that
$\forall i\in \N$,  $\E \|X_i \|^{2+\delta}<\infty$ and $\forall
k\in \N^*$, define
\begin{eqnarray}\label{condlinddelta}
A_k=\sum_{i=1}^k \E \|X_i\|^{2+\delta} .
\end{eqnarray}
Let $(Y_i)_{i\in \N}$ be a sequence of zero mean independent r.v.
with values in $\R^d$, independent of the sequence $(X_i)_{i\in
\N}$ and such that $Y_i\sim{\cal N}_d(0,\cov X_i)$ for all $i \in
\N$. Denote by ${\cal C}^3_b$ the set of bounded functions
$\R^d\to\R$ with bounded and continuous partial derivatives up to
order $3$. Set, for $f\in {\cal C}^3_b$ and $k\in \N^*$,
\begin{eqnarray}\label{methlind}
\Delta_k&=&\Big|\E\big(f(X_1+\cdots+X_k)-f(Y_1+\cdots+Y_k)\big)\Big|
\end{eqnarray}
Following the dependence between vectors $X_i$, we now provide $3$
lemmata, the first one is well known and relates to the
independence case, the two others are concerned with the
dependence case. Thus, first, for independent random variables,
the Lindeberg lemma (see e.g. Petrov, 1995) is
\begin{lemma}[Lindeberg under independence]\label{lemlindind}
Let $(X_i)_{i\in \N}$ be a sequence of independent zero mean r.v.
with values in $\R^d$ satisfying Assumption $H_\delta$. Then, for
all $k\in \N^*$:
$$
\Delta_k \le 3 \cdot \|f^{(2)}\|_\infty^{1-\delta} \cdot
\|f^{(3)}\|_\infty^\delta \cdot  A_k.
$$
\end{lemma}~\\
This lemma is restated for completeness sake but it is essentially
well known.
\begin{remark}\label{classic}
Using the proof of the previous Lemma, classical Lindeberg
conditions may be used:
\begin{equation}\label{lindclas}
\Delta_k \le \|f^{(2)}\|_\infty B_k(\varepsilon)+ \|f^{(3)}\|_\infty
\cdot a_k\left(\frac 2 3 \, \varepsilon+\frac 1 2
\,\sqrt{B_k(\varepsilon)}\right),
\end{equation}
\vspace{-0.5cm}
\begin{eqnarray}\label{condlind}
\mbox{where}~~~~B_k(\varepsilon)&=&\sum_{i=1}^k \E \Big (
\|X_i\|^2\1_{\{\|X_i\|>\varepsilon\}}\Big ),\ \mbox{ for
}\varepsilon>0,~k\in \N,
\\
\nonumber
 a_k&=&\sum_{i=1}^k \E (\|X_i\|^2)<\infty,\ \mbox{ for
}~k\in \N.
\end{eqnarray}
Moreover, these classical Lindeberg conditions derive those from
Lemma \ref{lemlindind}; indeed, according to H\"older Inequality,
$\E \Big ( \|X_i\|^2\1_{\{\|X_i\|>\varepsilon\}}\Big ) \leq \Big
(\E ( \|X_i\|^{2+\delta})\Big )^{\frac 2{2+\delta}} \Big ( \P (\|
X_i \|
> \varepsilon \Big )^{\frac \delta {2+\delta}}$, and then,
using the Bienaym\'e-Tchebichev inequality, for all $\delta \in
]0,1[$, $B_k(\varepsilon)\le \varepsilon^{-\delta} A_k$.
Consequently (\ref{lindclas}) implies,
$$
\Delta_k \le \|f^{(2)}\|_\infty \varepsilon^{-\delta} A_k +
\|f^{(3)}\|_\infty a_k\big(\frac 2 3 \, \varepsilon+\frac 1 2
\,\varepsilon^{-\delta/2} \sqrt{A_k}\big).
$$
\end{remark}
~\bigskip For the dependent case, the Lindeberg method provides
the two following lemmata. First, for random vectors
$W=(W_1,\ldots,W_d)'$ and $X=(X_1,\ldots,X_d)'\in \R^d$, we will
use the notations,
\begin{eqnarray*}
\left \{ \begin{array}{lcl} \cov
(f^{(1)}(W),X)&=&\displaystyle{\sum_{\ell=1}^d\cov
\left(\frac{\partial f}{\partial x_\ell}(W),X_\ell\right),}
\\
\nonumber \cov
(f^{(2)}(W),X^2)&=&\displaystyle{\sum_{k=1}^d\sum_{\ell=1}^d\cov
\left(\frac{\partial^2 f}{\partial x_k\partial
x_\ell}(W),X_kX_\ell\right).}
\end{array} \right.
\end{eqnarray*}
Assume that $X^*$ is independent of $W$ and admits the same
distribution as $X$; the previous expression is rewritten $ \cov
(f^{(2)}(W),X^2)=\E f^{(2)}(W)(X, X)-\E f^{(2)}(W)(X^*, X^*)$
where $f^{(2)}(W)$ is here considered as a quadratic form of
$\R^d$.  Then,
\begin{lemma}[Dependent Lindeberg Lemma - I]\label{lemlind}
Let $(X_i)_{i\in \N}$ be a sequence of zero mean r.v. with values in
$\R^d$ satisfying Assumption $H_\delta$. Then,
$$
\Delta_k \le T_1(k)+\frac12 \, T_2(k)+ 6 \cdot
\|f^{(2)}\|_\infty^{1-\delta} \cdot \|f^{(3)}\|_\infty^\delta \cdot A_k,
$$
where (empty sums are set equal to $0$),
\begin{eqnarray}
\label{somcovlind} T_j(k)&=&\sum_{i=1}^k \Big |
\cov\left(f^{(j)}(X_1+\cdots+X_{i-1}),X_i^{j}\right)\Big |, \qquad
j=1,2.
\end{eqnarray}
\end{lemma}
Characteristic functions are considered below; they provide a
simpler result
\begin{lemma}[Dependent Lindeberg Lemma - II]\label{lemlind'}
Let $(X_i)_{i\in \N}$ be a sequence of zero mean r.v. with values in
$\R^d$ satisfying Assumption $H_\delta$. For the special case of
complex exponential functions $f(x)=e^{i<t,x>}$ for some $t\in\R^d$
(and where $<a , b>$ is the scalar product in $\R^d$),
$$
\Delta_k \le T(k)+3 \|t\|^{2+\delta}
A_k,~~~\mbox{where}~~T(k)=\sum_{j=1}^k \big |
\cov(e^{i<t,X_1+\cdots+X_{j-1}>},e^{i<t,X_j>})\big |.
$$
\end{lemma}
The main consequence of those three lemmata is related to the
asymptotic behavior of $\sum_{i=1}^k X_i$, and provide sufficient
conditions for establishing the C.L.T.
\begin{theorem}[a Lindeberg C.L.T.]\label{theolind}
Assume that the sequence $(X_{i,k})_{i\in \N}$ satisfies
Assumption $H_\delta$, and $A_k \limitek 0$, and there exists
$\Sigma$ a positive matrix such that
$\displaystyle{\Sigma_k=\sum_{i=1}^k\cov (X_{i,k})\limitek \Sigma}$.
Moreover, assume that $T_j(k) \limitek 0$ for $j=1,2$ (Lemma
\ref{lemlind}) or $T(k) \limitek 0$ (Lemma \ref{lemlind'}). Then,
$$
S_k=\sum_{i=1}^k X_{i,k} \limiteloik {\cal N}_d(0,\Sigma).
$$
\end{theorem}
\begin{proof}[Proof of Theorem \ref{theolind}] Under the assumptions
of this Theorem, it is clear that
$\Big|\E\big(f(S_k)-f(N_k)\big)\Big| \limitek 0$ for all functions
$f\in {\cal C}^3_b$, or for all $t\in \R$ and $f(x)=e^{itx}$,
where $N_k \sim {\cal N}_d(0,\Sigma_k)$. According to
$\Big|\E\big(f(N_k)-f(N)\big)\Big| \limitek 0$ where $N \sim {\cal
N}_d(0,\Sigma)$, we deduce that $\Big|\E\big(f(S_k)-f(N)\big)\Big|
\limitek 0$ and therefore $S_k \limiteloik N$. $\qquad\square$
\end{proof}
Following this theorem, we can remark that the condition
$A(k)\limitek 0$ is the usual Lindeberg condition (with also
condition $\Sigma_k \limitek \Sigma$, the convergence of variances),
while the conditions $T_j(k)\limitek 0$ (for $j=1,2$) or
$T(k)\limitek 0$ are related to the  dependence structure of the
sequence $(X_i)_{i\in \N}$.
\section{Weakly dependent processes}\label{sec3}
We have just seen that the convergence in distribution of $S_k$ to
a Gaussian law is obtained if $T_1(k)$ and $T_2(k)$, or $T(k)$
converge to $0$. Those terms are related to the dependence of the
sequence $(X_n)_{n\in \N}$. Now, we address a very general class
of dependent processes introduced and developped in Doukhan and
Louhichi (1999). Numerous reasons justify this choice. First, this
frame of dependence includes a lot of models like causal or non
causal linear, bilinear, strong mixing processes or also dynamical
systems. Secondly, these properties of dependence are independent
of the marginal distribution of the time series, that can be as
well a discrete one, Lebesgue measurable one or else. Finally,
these definitions of dependence can be easily used in various
statistic contexts, in particular in the case of the establishment
of central limit theorems, since the previous bounds provided for
$T_1(k)$, $T_2(k)$ or $T(k)$ are
written with sums of covariances (see above).\\
~\\
To define a such dependent processes, first, for
$h:\left(\R^d\right)^u\to\R$ an arbitrary function, with $d,u \in
\N^*$, denote,
$$
\Lip h=\sup_{(y_1,\ldots,y_u)\ne(x_1,\ldots,x_u) \in (\R^d)^u }
\frac{\left| h(y_1,\ldots,y_u)- h(x_1,\ldots,x_u) \right|}
{\|y_1-x_1\|+\cdots+\|y_u-x_u\|}.$$
Then,
\begin{definition}
A process $X=(X_n)_{n\in \Z}$ with values in $\R^d$ is a so-called
$(\varepsilon,\psi)$-weakly dependent process if there exist a
function $\psi:(\N^*)^2\times(\R^+)^2\to\R^+$ and a sequence
$(\varepsilon_r)_{r\in \N}$ such that $\varepsilon_r\limiter 0$
satisfying,
\begin{equation} \label{eta_ine}
\left|{{\cov}}\Big( g_1(X_{i_1},\ldots,X_{i_u}),
g_2(X_{j_1},\ldots,X_{j_v})\Big ) \right| \leq \psi(u,v,\lip
g_1,\lip g_2)\cdot \varepsilon_{r}
\end{equation}
for all~$\displaystyle{\left \{ \begin{array}{l} \bullet~(u,v)\in
\N^*\times \N^*;\\
\bullet~(i_1,\ldots,i_u)\in \Z^u~\mbox{and}~(j_1,\ldots,j_v)\in
\Z^v~\mbox{with}~ i_1\leq \cdots \leq i_u < i_u+r\leq j_1\leq\cdots\leq
j_v \\
\bullet~\mbox{functions}~ g_1:\R^{ud} \to
\R~\mbox{and}~g_2:\R^{vd} \to
\R~\mbox{such that}\\
\hspace{5cm}\|g_1\|_\infty\le 1,~\|g_2\|_\infty\le 1,~\Lip\ g_1 < \infty~\mbox{and}~\Lip\ g_2 < \infty; \\
 \end{array} \right . }$
\end{definition}
In the sequel, two different particular cases of functions $\psi$
corresponding to two different cases of weakly dependent processes
will be considered (more details can be found in Doukhan and
Louhichi, 1999, Doukhan and Wintenberger, 2006),
\begin{itemize}
\item If $X$ is a causal time series, {\it i.e.} there exist a sequence
of functions $(F_n)$ and a sequence of independent random
variables $(\xi_k)_{k\in \Z}$ such that
$X_n=F_n(\xi_n,\xi_{n-1},\ldots)$ for $n \in \Z$, the
$\theta$-weakly dependent causal condition, for which
$$
\psi(u,v,\lip\ g_1,\lip\ g_2)= v \cdot \lip\ g_2
$$
(in such a case, we will simply denote $\theta_r$ instead of
$\varepsilon_r$).
\item If $X$ is a non causal time series, the $\lambda$-weakly dependent condition, for which
$$
\psi(u,v,\lip\ g_1,\lip\ g_2)=u\cdot v\cdot \lip\ g_1 \cdot \lip\
g_2+u\cdot \lip\ g_1+v \cdot \lip\ g_2
$$
(in such a case, we will simply denote $\lambda_r$ instead of
$\varepsilon_r$).
\end{itemize}
\begin{remark}
It is clear that if $X$ is a $\theta$-weakly dependent process it
is also a $\lambda$-weakly dependent process. The main reasons for
considering a distinction between causal and non causal time
series are: a/ the $\theta$-weak dependence is more easily relied
to the strong mixing property; b/ some models or properties
require different conditions on the convergence rate of
$(\theta_r)$ than for $(\lambda_r)$.
\end{remark}
Note first that sums of  independent weakly dependent processes
admit the common weak dependence property where dependence
coefficients are the sums of the initial ones. We now provide a non
exhaustive list of weakly dependent sequences with their weak
dependence properties. In the sequel, $X=(X_k)_{k\in \Z}$ denote a
weakly dependent stationary time series (the conditions of the
stationarity will not be specified) and $(\xi_n)_{n\in \Z}$ is a
sequence of zero mean i.i.d. random variables,
\begin{enumerate}
\item If $X$ is a Gaussian process and if
$\displaystyle{\lim _{i\to \infty}  \cov(X_0,X_i)=0}$, then $X$ is
a $\lambda$-weakly dependent process such that
$\displaystyle{\lambda_r=O\Big( \sup _{i \geq
r}|\cov(X_0,X_i)|\Big)}$ (see Doukhan and Louhichi, 1999).
\item If $X$ is an associated stationary processes, then $X$ is $\lambda$-weakly
dependent process such that $\displaystyle{\lambda_r=O\Big(
\sup_{i \geq r}\cov(X_0,X_i)\Big)}$ (see Doukhan and Louhichi,
1999).
\item if $X$ is a $ARMA(p,q)$ process or,
more generally, a causal (respectively, a non causal) linear
process such that $\displaystyle{X_k=\sum_{j=0}^\infty a_j
\xi_{k-j}}$ (respectively,
$\displaystyle{X_k=\sum_{j=-\infty}^\infty a_j \xi_{k-j}}$) for
$k\in \Z$, with $a_k={\cal O}(|k|^{-\mu})$ with $\displaystyle{\mu
>1/2}$, then $X$ is a $\theta$- (respectively, $\lambda$-) weakly dependent process with
$\displaystyle{\theta_{r}= \lambda_{r}={\cal O}\big(\frac 1
{r^{\mu-1/2}} \big )}$ (see Doukhan and Lang, 2002, p. 3). It is
also possible to deduce $\lambda$-weak dependence properties for
$X$ if the innovation process is itself $\lambda$-weakly dependent
(Doukhan and Wintenberger, 2006).
\item if $X$ is a $GARCH(p,q)$ process or, more generally, a $ARCH(\infty)$
process such that $\displaystyle{X_k= \rho_k \cdot \xi_k}$ with
 $\displaystyle{\rho_k^2 = b_0 + \sum_{j=1}^\infty
b_j X^2_{k-j}}$ for $k\in \Z$ and if,
\begin{itemize}
\item it exists $C>0$  and $\mu \in ]0,1[$ such that $\forall j \in \N$, $0 \leq b_j
\leq C \cdot \mu^{-j}$, then $X$ is a $\lambda$-weakly dependent
process with $\lambda_r={\cal O}(e^{-c \sqrt{r}})$ and $c>0$ (this
is the case of $GARCH(p,q)$ processes).
\item it exists $C>0$ and $\displaystyle{\nu >1}$
such that $\forall j \in \N$, $0 \leq b_j \leq C \cdot j^{-\nu}$, then
$X$ is a $\lambda$-weakly dependent process with
$\displaystyle{\lambda_r={\cal O}\big ( r^{-\nu+1} \big )}$ (see
Doukhan {\it et al}, 2005).
\end{itemize}
\item if $X$ is a causal bilinear process
such that $\displaystyle{X_k=\xi_k \Big ( a_0 + \sum_{j=1}^\infty
a_j X_{k-j} \Big ) +c_0+ \sum_{j=1}^\infty c_j X_{k-j}}$ for $k\in
\Z$ (see Giraitis and Surgailis, 2002) and if,
\begin{itemize}
\item $\displaystyle{\left \{ \begin{array}{l} \mbox{$\exists J \in \N$ such that $\forall j>J$,
$a_j=c_j=0$, or,}\\
\exists \mu \in ]0,1[~\mbox{such that $\sum _j |c_j| \mu ^{-j}
\leq 1$ and $\forall j \in \N$, $0 \leq a_j \leq \mu^j$}
\end{array} \right . ,}$ then $X$ is a $\lambda$-weakly dependent process with $\lambda_r={\cal
O}(e^{-c \sqrt{r}})$, constant $c>0$;
\item $\forall j \in \N$,
$c_j\geq 0$, and $\displaystyle{\exists \nu_1 >2}$ and $\exists
\nu_2>0$ such that $a_j={\cal O}(j^{-\nu_1})$ and $\sum _j c_j
j^{1+\nu_2}<\infty$, with $\displaystyle{d=\max \Big (
-(\nu_1-1)\,;\, -(\nu_2  \delta)(\delta+\nu_2  \log 2)^{-1}\Big
)}$, then $X$ is a $\lambda$-weakly dependent process with
$\displaystyle{\lambda_r={\cal O}\Big ( \big ( \frac r {\log r}
\big )^d \Big )}$ and  (see Doukhan {\it et al.}, 2004).
\end{itemize}
\item if $X$ is a non causal bilinear process satisfying
$\displaystyle{X_k=\xi_k \cdot \Big (a_0 + \sum_{j \in \Z^*} a_j
X_{k-j} \Big )}$, for $k \in \Z$, where
$\|\xi_0\|_{\infty}<\infty$ (bounded random variables) and
$a_k={\cal O}(|k|^{-\mu})$ with $\displaystyle{\mu>1}$, then $X$
is a $\lambda$-weakly dependent process with
$\displaystyle{\lambda_{r}= {\cal O}\big(r^{1-\mu} \big )}$ (See
Doukhan {\it et al.}, 2005).
\item  if $X$ is a non causal finite order Volterra process such
that $X_k=\sum_{p=1}^\infty Y_k^{(p)}$ for $k\in \Z$, and with
$\displaystyle{Y_k^{(p)}=\sum_{-\infty<j_1<j_2<\cdots<j_p<\infty}
a^{(p)}_{j_1,\ldots, j_p}\xi_{k-j_1}\cdots\xi_{k-j_p}}$ and such that
it exists $p_0\in \N^*$ satisfying for $p>p_0$,
$a^{(p)}_{j_1,\ldots,j_p}=0$. Then if
$\displaystyle{a^{(p)}_{j_1,\ldots,j_p}={\cal O} \Big (\max_{1\le
i\le p} \{ |j_i|^{-\mu} \} \Big )}$ with $\displaystyle{\mu >0}$,
$X$ is a $\lambda$-weakly dependent process with
$\displaystyle{\lambda_{r}= {\cal O}\big(\frac 1 {r^{\mu+1}} \big
)}$ (See Doukhan, 2003). As in case 3, $\lambda$-weak dependence
properties for $X$ may be proved even for $\lambda$-weakly
dependent innovations.
\item if $X$ is a causal (respectively, non causal) infinite memory
process such that
$$
X_k=F(X_{k-1},X_{k-2},\ldots;\xi_k)~~~~\mbox{(respectively,}~
X_k=F(X_{k-t},\,t\ne 0;\xi_k))~~~\mbox{for}~ k \in \Z,
$$
where the function $F$ is defined on $\R^{\N}$ (respectively,
$\R^{\Z}$) and satisfies, with $m >0$, $\|F(0;\xi_0)\|_m<\infty$
and $\|F((x_j)_j;\xi_0)-F((y_j)_j;\xi_0)\|_m\le
\sum_{j\ne0}a_j|x_j-y_j|$, where $a=\sum_{j\ne0}a_j<1$. Then $X$
is a $\theta$- (respectively, $\lambda$-) weakly dependent process
with $\theta_r=\inf_{p\ge1}\Big\{ a^{\frac
rp}+\sum_{j>p}a_j\Big\}$ (respectively,
$\lambda_r=\inf_{p\ge1}\Big\{ a^{\frac rp}+\sum_{|j|>p}a_j\Big\}$)
(see Doukhan and Wintenberger, 2006).
\item let $X$ be a  zero-mean, second order stationary Gaussian (or linear) process. Assume that $X$ is long-range
dependent:  $\cov(X_0,X_k)=L(k)\cdot k^{2H-2}$ for $k\in \N$, with
$H\in (1/2,1)$ (so-called Hurst parameter) and $L(\cdot)$ a slowly
varying function (at $+\infty$). Then, $X$ is a $\lambda$-weakly
dependent process with $\lambda_r=L(r)\cdot r^{2H-2}$.
\end{enumerate}
\noindent Now, different applications of Theorem \ref{theolind}
are considered for weakly dependent processes.
\section{Applications  under weak dependence}\label{appli}
\subsection{Lindeberg Central Limit Theorem}
Doukhan and Wintenberger (2006) prove a C.L.T. using Bernstein
blocks for a sequence $(X_i)_{i\in \N}$ of stationary zero mean
$(2+\delta)$-order random variables. In order to prove
$T(k)\limitek 0$ (see Lemma \ref{lemlind'}), we consider two
sequences $(p_n)_{n \in \N}$ and $(q_n)_{n \in \N}$ such that,
$$
\left \{ \begin{array}{l} p_n \limitek \infty \\
q_n \limitek \infty
\end{array} \right .~~\mbox{and}~~p_n=o(n),~~q_n=o(p_n).
$$
Introduce now the number of Bernstein blocks
$\displaystyle{k_n=\left[\frac n{p_n+q_n}\right]}$. It can be
shown that if $p_n\cdot q_n=o(n)$,
$$
\left \| \frac1{\sqrt n}\sum_{i=1}^nX_i-\frac1{\sqrt n}
\sum_{j=1}^{k_n} \sum_{i=(j-1)(p_n+q_n)+1}^{(j-1)(p_n+q_n)+p_n}
X_i \right \|_2 \limiten 0.
$$
Therefore, it is sufficient to prove the convergence in
distribution to a Gaussian law of the second sum, which is easier
than with the first one. Thus, Doukhan and Wintenberger (2006)
prove a $(2+\delta)$-order moment inequality which entails
condition $A(k_n)\to 0$ and $T(k_n)\limitek 0$, and they obtain:
\begin{theorem}
Let $(X_i)_{i\in \N}$ be a sequence of stationary zero mean
$(2+\delta)$-order random variables (with $\delta>0$). Assume that
$(X_i)_{i\in \N}$ is a $\lambda$- (or $\theta$-) weakly dependent
time series satisfying $\lambda_r=O(r^{-c})$ (or
$\theta_r=O(r^{-c})$) when $r \to \infty$, with $c >4+2/\delta$,
then it exists $0<\sigma^2<\infty$ such that
$\displaystyle{\frac1{\sqrt n}\sum_{i=1}^nX_i \limiteloik {\cal
N}(0,\sigma^2)}$.
\end{theorem}
Note that in Doukhan and Wintenberger (2006) the Donsker principle
is also proved.
\subsection{Subsampling}
Assume that $(X_i)_{i\in\Z}$ is a zero mean $(2+\delta)$-order
stationary sequence for some $\delta>0$, with $\Sigma=\cov (X_0)$.
Then, for a sequence $(m_n)_{n\in \N}$ such that $m_n \limiten
\infty$ and $\displaystyle{k_n=\big[n /{m_n}\big]\limiten \infty}$.
We consider a subsample $(X_{m_n},\ldots,X_{k_nm_n})$ of
$(X_1,\ldots,X_n)$, and the sample,
$$
(Y_{1,n},\ldots,Y_{k_n,n})~~\mbox{with}~~Y_{i,n}=\frac 1 {\sqrt{
k_n}} \,  X_{im_n}~~\mbox{for $1\le i \le k_n$.}
$$
Depending on the weak dependence property of $(X_i)_{i\in\Z}$, we
can obtain the Lindeberg Theorem for
$$
S_{k_n,n}=\sum_{i=1}^{k_n} Y_{i,n}.
$$
\begin{prop}\label{subsampling}
Assume that $(X_i)_{i\in\Z}$ is a zero mean $(2+\delta)$-order
stationary sequence for some $\delta>0$, with $\Sigma=\cov (X_0)$.
Then, for a sequence $(m_n)_{n\in \N}$ such that $m_n \limiten
\infty$ and $\displaystyle{k_n=\big[n/ {m_n}\big]\limiten
\infty}$,
$$
S_{k_n,n}=\frac1{\sqrt{ k_n}}\,\sum_{i=1}^{k_n}X_{im_n}\limiteloin
{\cal N}_d(0,\Sigma),
$$
if one of the following assumptions also holds,
\begin{itemize}
\item $(X_i)_{i\in\Z}$ is a $\theta-$weakly dependent sequence and
$\theta_{m_n}\sqrt{ k_n} \limiten 0$.
\item $(X_i)_{i\in\Z}$ is a $\lambda-$weakly dependent sequence and
$\lambda_{m_n}k_n^{\frac32}\limiten 0$.
\end{itemize}
\end{prop}
This proposition allows to pass from a situation where the C.L.T.
is not satisfied to a situation where it is satisfied by using a
subsample with the correct asymptotic step of sampling. For
instance, if $X$ is a zero mean stationary long range dependent
$(2+\delta)$-order process (with $\delta>0$) and such that a/ $X$
is a Gaussian process such that $\E (X_0X_n)=O(n ^{2H-2})$ with
$1/2<H<1$ when $n \to \infty$ or b/ $X$ a linear process. Then $X$
is a $\lambda$-weakly dependent process with
$\lambda_r=O(r^{2H-2})$ and it is well known (see for instance
Taqqu, 1975) that $X$ does not satisfy a usual central limit
theorem. As a consequence, with a subsampling step $m_n$ such that
$o(m_n)=n^{3/(4H-1)}$, the subsampled time series $(X_{jm_n})$
satisfies a usual C.L.T. with a convergence rate $o(n^{(1-H)/(4H-1)})$. \\
~\\
Two objections can be raised to this method: first, only a part of
the sample is used. The second and main objection is that the
choice of the convergence rate of the subsampling implies the
knowledge of the convergence rate of $(\lambda_r)$ or
$(\theta_r)$. However, an estimation of this last rate (for
instance in the case of long-range dependence) could provide an
estimation of a fitted step of subsampling, or in the case of an
exponential rate of convergence of $(\lambda_r)$ or $(\theta_r)$,
all subsampling step $m_n=O(n^a)$ with $0<a<1$ could be used: it
is the case for instance for $GARCH(p,q)$ processes.
\subsection{Kernel density estimation}
Let $(X_i)_{i\in \N}$ be a sequence of stationary zero mean r. v.
(with real values) such that $X_0$ has a marginal density $f_X$
with respect to Lebesgue measure. Let $K:\R\to\R$ be
a kernel function satisfying,\\
~\\
{\bf Assumption K:} $K:\R\to\R$ be a bounded and Lipschitz
function with $\int_{-\infty}^\infty K(t)\,dt=1$.\\
~\\
Now, define $(h_n)_{n\in \N}$ a sequence such that $h_n \limiten
0$ and consider the usual kernel density estimation,
$$
\widehat f_X^{(n)}(x)=\frac 1 n \sum_{i=1}^n \frac 1 {h_n} \,
K\left(\frac{x-X_i}{h_n}\right)~~~\mbox{for $x \in \R$}.
$$
\begin{prop}\label{kernel}
Let $(X_i)_{i\in \Z}$ be a stationary zero mean weakly dependent
time series (with real values) such that $X_0$ has a marginal
density $f_X$ with respect to Lebesgue measure. Assume that
$\|f_X\|<\infty$ and $\max_{i\ne j}\|f_{i,j}\|_\infty<\infty$,
where $f_{i,j}$ denotes the joint marginal density of $(X_i,X_j)$.
Then,
$$
\sqrt{nh_n}\, \Big (\widehat f_X^{(n)}(x)-\E \widehat f_X^{(n)}(x)
\Big ) \limiteloin {\cal N}\Big (0,f_X(x)\int_\R K^2(t)\,dt\Big ),
$$
if $h_n \limiten 0$, $nh_n \limiten \infty$ and,
\begin{itemize}
\item $h_n=o\big (n^{-2/(\lambda-4)}\wedge n^{-5/(2\lambda-5) }\big )$ when $(X_i)_{i\in \Z}$ is a
$\lambda$-weakly dependent process with
$\lambda_r=O(r^{-\lambda})$ and $\lambda>5$, or
\item when $(X_i)_{i\in \Z}$ is a
$\theta$-weakly dependent process with $\theta_r=O(r^{-\theta})$
and $\theta>3$.
\end{itemize}
\end{prop}
Now, following the regularity of the function $f_X$, $\E \widehat
f_X^{(n)}(x)$ is a more or less good approximation of $f_X(x)$.
Hence, here there are two different cases of the approximation of
the bias,
\begin{cor}\label{regu1}
Assume that the function $f_X$ is a ${\cal C}^p(\R)$ function,
with $p \in \N^*$. Then, under the same conditions than
Proposition \ref{kernel}, with $K$ a kernel such that $\int _\R
K(t)t^q\,dt=0$ for $q=\{1,\ldots, p-1\}$ and $\int _\R K(t)t^p\,dt
\neq 0$, if $h_n=C \cdot n^{-1/(2p+1)}$ (with $C>0$) and, following
the different frames of weak dependence, $\lambda >5p+5$ or
$\theta>3$,
$$
\sqrt{nh_n}\, \Big (\widehat f_X^{(n)}(x)-f_X(x) \Big )
\limiteloin {\cal N}\Big (f_X^{(p)}(x)\frac 1 {p !} \,  \int _\R
t^{p} K(t)\, dt  \, ,\, f_X(x)\int_\R K^2(t)\,dt\Big ).
$$
\end{cor}
\begin{cor}\label{regu2}
Assume that the regularity of the function $f_X$ is $\rho>0$ with
$\rho \notin \N^*$ (in the sense that $f$ has a $[\rho]$-th order
bounded derivative which is H\"older continuous with exponent
$\rho-[\rho]$). Then, under the same conditions than Proposition
\ref{kernel}, with $K$ a kernel such that $\int _\R K(t)t^q\,dt=0$
for $q=\{1,\ldots, [\rho]-1\}$ and $\int _\R K(t)t^{[\rho]}\,dt
\neq 0$, if $h_n=o\big (n^{-1/(2[\rho]+1)}\big )$ and, following
the different frames of weak dependence, $\lambda >5[\rho]+5$ or
$\theta>3$,
$$
\sqrt{nh_n}\, \Big (\widehat f_X^{(n)}(x)-f_X(x)\Big ) \limiteloin
{\cal N}\Big (0,f_X(x)\int_\R K^2(t)\,dt\Big ).
$$
\end{cor}
\subsection{Subsampled kernel density estimation}
Now, imagine that the process $(X_k)_{k\in \Z}$ is a weakly
dependent zero mean stationary process such that the conditions
$\lambda> 5$ or $\theta > 3$ of Proposition \ref{kernel} are not
satisfied. As a consequence, the kernel density estimator is not
proved to satisfy a central limit theorem. Subsampling can provide
a way for obtaining a C.L.T. (and then, confidence intervals or
goodness-of-fit tests). Indeed, like it was also considered
before, a subsampled time series with an asymptotic rate is "less"
dependent than the original time series. Indeed, consider a
sequence $(m_n)_{n\in \N}$ such that
$$
m_n \limiten \infty~~~\mbox{and}~~~k_n=\big[n /{m_n}\big]\limiten
\infty,
$$
and the subsample $(X_{m_n},\ldots,X_{k_nm_n})$ of
$(X_1,\ldots,X_n)$. For $(h_n)_{n\in \N}$ a sequence such that
$h_n \limiten 0$, define the subsampled kernel density estimator
of $f_X$ as:
$$
\widehat f_X^{(n,m_n)}(x)=\frac 1 {k_n} \sum_{i=1}^{k_n} \frac 1
{h_n} \, K\left(\frac{x-X_{im_n}}{h_n}\right)~~~\mbox{for $x \in
\R$}.
$$
\begin{prop}\label{subkernel}
Under the same assumptions than Proposition \ref{kernel}, except
that $0\leq \lambda \leq 6$ or $\theta\leq 3$, the following C.L.T.
yields from the subsample $(X_{m_n},\ldots,X_{k_nm_n})$,
$$
\sqrt{k_n h_n}\, \Big (\widehat f_X^{(n,m_n)}(x)-\E \widehat
f_X^{(n,m_n)}(x) \Big ) \limiteloin {\cal N}\Big (0,f_X(x)\int_\R
K^2(t)\,dt\Big ),
$$
for sequences $(h_n)$ and $(m_n)$ such that $h_n=n^{-h}$ and
$m_n=n^m$, where the points $(m,h)$ are in the "white zonas"
(included in square $(0,1)^2$) of Figure 1 and 2 (see below).
Moreover, the "optimal" convergence rate is obtained for,
\begin{itemize}
\item $\displaystyle{\sqrt{h_nk_n}=n^{\frac \lambda {5+2(\lambda \vee 1)}-\varepsilon}}$
for all $\varepsilon >0$ small enough (and $h_n=n^{-\varepsilon}$
and $\displaystyle{m_n=n^{\frac 5 {5+2(\lambda \vee
1)}+\varepsilon}}$) when $(X_i)_{i\in \Z}$ is a $\lambda$-weakly
dependent process;
\item $\displaystyle{\sqrt{h_nk_n}=n^{\frac 1 2 (\theta \wedge 1)-\varepsilon}}$
for all $\varepsilon >0$ small enough (and $h_n=n^{-\varepsilon}$
and $\displaystyle{m_n=n^{((1-\theta)\vee 0)+\varepsilon}}$) when
$(X_i)_{i\in \Z}$ is a $\theta$-weakly dependent process.
\end{itemize}
\end{prop}
From this result, we now answer to the following problem. Assume
that $X$ is a $\theta$- (or $\lambda$-) weakly dependent process and
that the regularity $\rho$ (in the sense of the previous section) of
the density function $f_X$ and the sequence decay rates $\theta$ and
$\lambda$ of the sequences $(\theta_r)_r$ or $(\lambda_r)_r$ are
known. In the previous section, we have established a C.L.T. for the
density kernel estimator when $\theta$ (or $\lambda$) is larger than
an affine function of $\rho$. However, if this inequality is not
satisfied, it is possible to find a subsampling step such that a
C.L.T. for the subsampled kernel density estimator can be proved,
\begin{cor}
Under the assumptions of consequence \ref{regu1}, but if
$0<\lambda \leq 5p+5$ (or $0< \theta \leq p+3$),
$$
\sqrt{k_nh_n}\, \Big (\widehat f_X^{(n,m_n)}(x)-f_X(x)\Big )
\limiteloin {\cal N}\Big (f_X^{(p)}(x)\frac 1 {p !} \,  \int _\R
t^{p} K(t)\, dt \,,\, f_X(x)\int_\R K^2(t)\,dt\Big ),
$$
with the following optimal conditions,
\begin{itemize}
\item for $\lambda$-weakly dependent process, with
$\displaystyle{\lambda_r=O(r^{-\lambda})}$ and $0<\lambda \leq
5p+5$, for $\displaystyle{h_n=n^{-\frac \lambda {5+p(5+2(\lambda
\vee 1))}}}$ and $\displaystyle{m_n=n^{1-\frac
{\lambda(2p+1)}{5+p(5+2(\lambda \vee 1))}}}$. Then, the
convergence rate is $\displaystyle{\sqrt{k_nh_n} \sim n^{\frac
{p\lambda}{5+p(5+2(\lambda \vee 1))}}}$;
\item for $\theta$-weakly dependent process, with
$\displaystyle{\theta_r=O(r^{-\theta})}$ and $0<\theta \leq p+3$,
for $\displaystyle{h_n=n^{-\frac \theta {3+2p(\theta\vee 1) }}}$
and $\displaystyle{m_n=n^{1-\frac {\theta(2p+1)}{3+2p(\theta\vee
1)}}}$. Then, the convergence rate is $\displaystyle{\sqrt{k_nh_n}
\sim n^{\frac {p\theta}{3+2p(\theta\vee 1)}}}$.
\end{itemize}
\end{cor}
\begin{cor}
Under the assumptions of Consequence \ref{regu2} with $\rho \notin
\N$, but if $0<\lambda \leq 5[\rho]+5$ (or $0< \theta \leq
[\rho]+3$), then,
$$
\sqrt{k_nh_n}\, \Big (\widehat f_X^{(n,m_n)}(x)-f_X(x)\Big )
\limiteloin {\cal N}\Big (0 \,,\, f_X(x)\int_\R K^2(t)\,dt\Big ),
$$
with the following conditions (for all $\varepsilon>0$ small
enough),
\begin{itemize}
\item for $\lambda$-weakly dependent process, with
$\displaystyle{\lambda_r=O(r^{-\lambda})}$ and $0<\lambda \leq
5[\rho]+5$, for $\displaystyle{h_n=n^{-\frac \lambda
{5+[\rho](5+2(\lambda \vee 1))}}}$ and
$\displaystyle{m_n=n^{1+2\varepsilon-\frac
{\lambda(2[\rho]+1)}{5+[\rho](5+2(\lambda \vee 1))}}}$. Then, the
convergence rate is $\displaystyle{\sqrt{k_nh_n} \sim n^{\frac
{p\lambda}{5+[\rho](5+2(\lambda \vee 1))}-\varepsilon}}$;
\item for $\theta$-weakly dependent process, with
$\displaystyle{\theta_r=O(r^{-\theta})}$ and $0<\theta \leq
[\rho]+3$, for $\displaystyle{h_n=n^{-\frac \theta
{3+2[\rho](\theta\vee 1) }}}$ and
$\displaystyle{m_n=n^{1+2\varepsilon-\frac
{\theta(2[\rho]+1)}{3+2[\rho](\theta\vee 1)}}}$. Then, the
convergence rate is $\displaystyle{\sqrt{k_nh_n} \sim n^{\frac
{[\rho]\theta}{3+2[\rho](\theta\vee 1)}-\varepsilon}}$.
\end{itemize}
\end{cor}
Hence, for all regularity parameter $\rho >0$, even if $\lambda$
or $\theta$ are very small numbers (for instance when $X$ is a
long range dependent process), the subsampled kernel density
estimator satisfies a C.L.T. for a fitted choice of $h_n$ and
$m_n$. Moreover, for $\theta$-weakly dependent time series with
$\theta>0$, when $\rho \to \infty$, the convergence rate of this
theorem is $n^{1/2-\varepsilon}$, with $\varepsilon>0$.
\section{Proofs}\label{proofs}
\begin{proof}[Proof of Lemma \ref{lemlindind}]
For $k \in \N^*$, we first
notice and prove that,
\begin{eqnarray}\label{decomplind}
\Delta_k&\le& \Delta_{k,1}+\cdots+\Delta_{k,k}\\
\nonumber \mbox{with}~~~~\Delta_{k,
i}&=&\Big|\E\big(f_i(W_i+X_i)-f_i(W_i+Y_i)\big)\Big|,
~~~~\mbox{for $i\in \{1,\ldots,k\}$}\\
\nonumber W_i&=&X_1+\cdots+X_{i-1}~~~\mbox{and}~~~W_1=0,~~\mbox{for
$i\in \{2,\ldots,k\}$},
\\
\nonumber f_i(t)&=&\E \big ( f(t+Y_{i+1}+\cdots+Y_k) \big )
~~~\mbox{and}~~~f_k(t)=f(t),~~\mbox{for $t\in \R^d$ and $i\in
\{1,\ldots,k-1\}$}.
\end{eqnarray}
Let $x,w\in\R^d$. Taylor formula writes in two following ways (for
suitable vectors $w_{1,x},w_{2,x}\in\R^d$),
\begin{eqnarray*}
f(w+x)&=&f(w)+f^{(1)}(w)(x)+\frac12 \, f^{(2)}(w_{1,x})(x,x)
\\
&=&f(w)+f^{(1)}(w)(x)+\frac12 \, f^{(2)}(w)(x,x)+\frac16\,
f^{(3)}(w_{2,x})(x,x,x),
\end{eqnarray*}
where, for $j=1,2$ and $3$, $f^{(j)}(w)(y_1,\ldots, y_j)$ stands for
the value of the symmetric $j$-linear form $f^{(j)}$ of
$(y_1,\ldots, y_j)$ at $w$. Moreover, denote,
\begin{eqnarray*}
\|f^{(j)}(w)\|_1=\sup_{\|y_1\|,\ldots,\|y_j\|\le1}|f^{(j)}(w)(y_1,\ldots,
y_j)|,\quad \|f^{(j)}\|_\infty=\sup_{w\in \R^d}\|f^{(j)}(w)\|_1.
\end{eqnarray*}
Thus for $w,x,y\in \R^d$, there exists some suitable vectors
$w_{1,x},w_{2,x},w_{1,y},w_{2,y}\in\R^d$ such that:
\begin{eqnarray*}
f(w+x)-f(w+y)&=&f^{(1)}(w)(x-y)+\frac12\, \Big(f^{(2)}(w)(x,x)-f^{(2)}(w)(y,y)\Big)\\
&&\qquad+\frac 1 2 \, \Big (
(f^{(2)}(w_{1,x})-f^{(2)}(w))(x,x)-(f^{(2)}(w_{1,y})-f^{(2)}(w))(y,y)
\Big )\\
&=&f^{(1)}(w)(x-y)+\frac12\, \Big(f^{(2)}(w)(x,x)-f^{(2)}(w)(y,y)\Big)\\
&&\qquad+\frac16 \, \Big (
f^{(3)}(w_{2,x})(x,x,x)-f^{(3)}(w_{2,y})(y,y,y) \Big ).
\end{eqnarray*}
$\mbox{Thus},~\gamma=f(w+x)-f(w+y)-f^{(1)}(w)(x-y)-\frac12 \,
\left(f^{(2)}(w)(x,x)-f^{(2)}(w)(y,y)\right)~~\mbox{ satisfies}$
\begin{eqnarray}
\nonumber |\gamma|&\le& \Big ((\|x\|^2+\|y\|^2)\|f^{(2)}\|_\infty
\Big ) \wedge \Big (\frac1 6 \,(\|x\|^3+\|y\|^3)\|f^{(3)}\|_\infty
\Big )
\\
\label{majogamma2} &\le&\|f^{(2)}\|_\infty
\left\{\|x\|^2\Big(1\wedge \frac 1 {6} \,
\frac{\|f^{(3)}\|_\infty}{\|f^{(2)}\|_\infty}\|x\|\Big)
+\|y\|^2\Big(1\wedge \frac 1 {6} \,
\frac{\|f^{(3)}\|_\infty}{\|f^{(2)}\|_\infty}\|y\|\Big)\right\}
\\
\label{majogamma} &\le& \frac 1 {6^\delta}\,
\|f^{(2)}\|_\infty^{1-\delta}
\|f^{(3)}\|_\infty^\delta\left\{\|x\|^{2+\delta}+\|y\|^{2+\delta}\right\}
\end{eqnarray}
using the inequality $1\wedge c\le c^\delta$, that is valid for
all $c\ge0$ when $\delta \in [0,1]$. Now, this inequality could be
applied to the functions $f_i$ and r.v. $Y_j$ and $W_k$. Indeed,
\begin{multline*}
f_i(W_i+X_i)-f_i(W_i+Y_i)-f^{(1)}(W_i)(X_i-Y_i)-\frac12 \,
\left(f^{(2)}(W_i)(X_i,X_i)-f^{(2)}(W_i)(Y_i,Y_i)\right)\\
=f_i(W_i+X_i)-f_i(W_i+Y_i)
\end{multline*}
because $W_i$ is independent from $X_i$ and $Y_i$ and because
$\E(X_i)=\E(Y_i)=0$ and $\cov(X_i)=\cov(Y_i)$. Thanks to the
inequalities $ \|f_i^{(j)}\|_\infty\le\|f^{(j)}\|_\infty $ (valid
for $1\le i\le k$ and $0\le j\le 3$) and (\ref{majogamma}), we
obtain,
$$
\Delta_{k,i} \leq \frac 1 {6^\delta}\,
\|f^{(2)}\|_\infty^{1-\delta}
\|f^{(3)}\|_\infty^\delta\left\{\E\|X_i\|^{2+\delta}+\E\|Y_i\|^{2+\delta}\right\}.
$$
But, $\E(\|Y_i\|^{2+\delta}) \leq (\E
(\|Y_i\|^4))^{\frac12+\frac\delta4}$ and $\E (\|Y_i\|^4)=3 \cdot
\E^2(\|X_i\|^2)$ because $Y_i$ is a Gaussian r.v. with the same
covariance than $X_i$. Therefore,  $\E(\|Y_i\|^{2+\delta}) \leq
3^{\frac12+\frac\delta4} (\E(\|X_i\|^{2+\delta}))^{4
(\frac12+\frac\delta4)/(2+\delta)}\leq 3^{\frac12+\frac\delta4}
\E(\|X_i\|^{2+\delta})$. As a consequence, from assumption
(\ref{condlinddelta}),
\begin{eqnarray*}
\Delta_k &\le& \frac 1 {6^\delta}\|f^{(2)}\|_\infty^{1-\delta}
\|f^{(3)}\|_\infty^\delta
\sum_{i=1}^k\left\{\E(\|X_i\|^{2+\delta})+\E
(\|Y_i\|^{2+\delta})\right\}
\\
&\le& \frac {(1+ 3^{\frac12+\frac\delta4})}{6^\delta}\, A_k \cdot
\|f^{(2)}\|_\infty^{1-\delta} \|f^{(3)}\|_\infty^\delta
\\
&\le& 3\cdot A_k \cdot \|f^{(2)}\|_\infty^{1-\delta}
\|f^{(3)}\|_\infty^\delta. \qquad\square
\end{eqnarray*}
\end{proof}
\begin{proof}[Proof of Remark \ref{classic}]\label{proofclassic}
Set $b_k^2=\max_{1\le
i\le k}\E(\|X_i\|^2)$. Now, for $\varepsilon<6\, \|f^{(2)}\|_\infty
\cdot (\|f^{(3)}\|_\infty)^{-1}$, and using inequality
(\ref{majogamma2}) and $1\wedge c\le c$ for all $c>0$, one obtains
\begin{eqnarray*}
\Delta_{k,i} &\le&   \|f^{(2)}\|_\infty \E \Big ( \|X_i\|^2 \wedge
\frac1{6} \,  \big(\|f^{(2)}\|_\infty \big ) ^{-1}
\|f^{(3)}\|_\infty \cdot  \|X_i\|^3 \Big )\Big)  +\frac1 6\,
\|f^{(3)}\|_\infty\E (\|Y_i\|^3) \\
&\le&   \|f^{(2)}\|_\infty \E \big ( \|X_i\|^2 \1_{\{\|X_i\|>
\varepsilon \}}\big ) + \frac1{6} \, \|f^{(3)}\|_\infty \cdot \Big (
\E \big ( \|X_i\|^3 \1_{\{\|X_i\| \leq  \varepsilon  \}}\big ) +\E
(\|Y_i\|^3) \Big ) \\
&\le&   \|f^{(2)}\|_\infty \E \big ( \|X_i\|^2 \1_{\{\|X_i\|>
\varepsilon \}}\big ) + \frac1{6} \, \|f^{(3)}\|_\infty \cdot \Big
(\varepsilon \cdot  \E ( \|X_i\|^2 ) +3 ^{3/4} (\E
((\|X_i\|^2)^{3/2}) \Big ),
\end{eqnarray*}
from the H\"older Inequality. It implies that,
\begin{eqnarray*}
\Delta_k &\le& \|f^{(2)}\|_\infty B_k(\varepsilon)+\frac1{6}\,
\|f^{(3)}\|_\infty \Big (  \varepsilon \cdot  a_k +3 ^{3/4}
\sum_{i=1}^k b_k \cdot \E (\|X_i\|^2)
\Big ) \\
&\le & \|f^{(2)}\|_\infty B_k(\varepsilon)+\frac 1 6 \,
\|f^{(3)}\|_\infty \cdot a_k\big(\varepsilon+3^{3/4}\cdot b_k
\big).
\end{eqnarray*}
Moreover, $b_k^2\le \varepsilon^2+\max_{1 \leq i\le k}\E \Big (
\|X_i\|^2\1_{\{\|X_i\|>\varepsilon\}}\Big )$, therefore $b_k^2
\leq \varepsilon^2+{B_k(\varepsilon)}$ and thus
$\displaystyle{b_k\le  \varepsilon+\sqrt{B_k(\varepsilon)}}$. As a
consequence,
\begin{eqnarray*}
\Delta_k &\le&  \|f^{(2)}\|_\infty B_k(\varepsilon)+
\|f^{(3)}\|_\infty \cdot a_k\Big(\frac 2 3 \, \varepsilon+\frac 1 2 \,
\sqrt{B_k(\varepsilon)}\Big). \qquad\square
\end{eqnarray*}
\end{proof}
\begin{proof}[Proof of Lemma \ref{lemlind}]
Consider $(X_i^*)_{i\in
\N}$ a sequence of r.v. satisfying Assumption $H_\delta$ and such
that $(X_i^*)_{i\in \N}$ is independent of $(X_i)_{i\in \N}$ and
$(Y_i)_{i\in \N}$. Moreover, assume that $X_i^*$ has the same
distribution as $X_i$ for $i\in \N$. Then, using the same
decomposition as in the proof of Lemma \ref{lemlindind}, one can
also write,
$$
\Delta_{k,i} \le
\left|\E\left(f_i(W_i+X_i)-f_i(W_i+X_i^*)\right)\right|
+\left|\E\left(f_i(W_i+X_i^*)-f_i(W_i+Y_i))\right)\right|.
$$
From the previous Lemma,
$$
\sum_{i=1}^k
\left|\E\left(f_i(W_i+X_i^*)-f_i(W_i+Y_i)\right)\right| \leq 3
\cdot \|f^{(2)}\|_\infty^{1-\delta} \cdot
\|f^{(3)}\|_\infty^\delta \cdot A_k.
$$
Moreover, the bound for $\gamma$ of the proof Lemma \ref{lemlindind}
entails,
\begin{multline*}
\left|f_i(W_i+X_i)-f_i(W_i+X_i^*)\right|\leq \Big
|f_i^{(1)}(W_i)(X_i)\Big |+ \frac12 \, \Big |
f_i^{(2)}(W_i)(X_i,X_i)-f_i^{(2)}(W_i)(X_i^*,X_i^*) \Big |+ \\
+\frac 1 {6^\delta}\, \|f^{(2)}\|_\infty^{1-\delta}
\|f^{(3)}\|_\infty^\delta\left\{\|X_i\|^{2+\delta}+\|X_i^*\|^{2+\delta}\right\},
\end{multline*}
because $(X_i)$ is now supposed to be a dependent sequence of
random variables and is no more independent from $(W_i)$. But with
the notation before Lemma \ref{lemlind}, we may write
$$
\sum_{i=1}^k |\E f_i^{(1)}(W_i)(X_i)\Big | = T_1~~\mbox{and}~~
\sum_{i=1}^k \Big
|\E\left(f_i^{(2)}(W_i)(X_i,X_i)-f_i^{(2)}(W_i)(X_i^*,X_i^*)\right)\Big
| = T_2.
$$
It implies that,
$$
\Delta_k \le T_1+\frac12 \, T_2+ 6 \cdot
\|f^{(2)}\|_\infty^{1-\delta} \cdot \|f^{(3)}\|_\infty^\delta
\cdot A_k. \qquad\square
$$
\end{proof}
\begin{proof}[Proof of Lemma \ref{lemlind'}]
Here, for $t\in \R^d$,
$f(x)=e^{i<t,x>}$. Denote $V_i=\v X_i$, the covariance matrix of
the vector $X_i$. Then, for a r.v. $Z$ independent from
$(Y_i)_{i\in \N}$,
$$
\E f_j(Z)=\E \big ( f(Z+Y_{j+1}+\cdots+Y_k) \big )=e^{-\frac12
t'\cdot \left(V_{j+1}+\cdots+V_k\right)\cdot t}\cdot \E
(e^{i<t,Z>}).
$$
Then, again
$$
\Delta_{k,j} \le \left|\E(f_j(W_j+X_j)-f_j(W_j+X_j^*))\right|
+\left|\E(f_j(W_j+X_j^*)-f_j(W_j+Y_j))\right|,
$$
with the second term bounded as in the proof of Lemma
\ref{lemlind} with $\|f^{(2)}\|_\infty^{1-\delta} \cdot
\|f^{(3)}\|_\infty^\delta=\|t\|^{2+\delta}$, and for the first
term,
\begin{eqnarray*}
\E(f_j(W_j+X_j)-f_j(W_j+X_j^*))&=& e^{-\frac12 t'\cdot
\left(V_{j+1}+\cdots+V_k\right)\cdot t}\cdot \E \Big
(e^{i<t,W_j>}(e^{i<t,X_j>}-e^{i<t,X_j^*>})\Big ) , \\
\big | \E(f_j(W_j+X_j)-f_j(W_j+X_j^*))\big | &\leq & \left|
e^{-\frac12 t'\cdot \left(V_{j+1}+\cdots+V_k\right)\cdot
t}\right|\big | \cov(e^{i<t,W_j>},e^{i<t,X_j>})\big
|.\\
&\leq &\big | \cov(e^{i<t,W_j>},e^{i<t,X_j>})\big |.
\qquad\square
\end{eqnarray*}
\end{proof}
\begin{proof}[Proof of Proposition \ref{subsampling}]
The proof of this proposition is a consequence of Theorem
\ref{theolind}, using $T(n) \limiten 0$. In one hand, thanks to
the stationarity of the sequence $(X_i)_{i\in \Z}$,
\begin{eqnarray*}
A_{k_n} = \sum_{i=1}^{k_n} \E \big (\|Y_{i,n}\|^{2+\delta} \big )=
k_n^{-\delta/ 2}\E \big ( \|X_0\|\big )^{2+\delta} \limiten 0.
\end{eqnarray*}
In the other hand, a bound of $T(k_n)$ can be provided. Indeed, let
$t\in \R^d$, and then,
\begin{eqnarray*}
T(k_n)&=&  \sum_{j=1}^{k_n-1} \left |\cov\left(e^{i <t,
Y_1+\cdots+Y_{j-1}>},e^{ i <t,Y_j>}\right) \right |,\\
&=& \sum_{j=1}^{k_n-1} \left |\cov(F_{t,n}(X_1,\ldots,
X_{j-1}),G_{t,n}(X_{j})) \right |,
\end{eqnarray*}
with $G_{t,n}(s_j)=\exp (i<s,t>/\sqrt{k_n})$ and
$F_{t,n}(s_1,\ldots,s_{j-1})=G_{t,n}(s_1) \times \cdots  \times
G_{t,n}(s_{j-1})$ for $(s_1,\ldots,s_j) \in (\R^d)^j$. But,
\begin{eqnarray*}
\|G_{t,n}\|_\infty\le1~~&\mbox{and}&~~\Lip G_{t,n} \le
\|t\|\cdot  k_n^{-1/2} \\
\|F_{t,n}\|_\infty \le1 ~~&\mbox{and}&~~\Lip F_{t,n} \le
\|t\|\cdot k_n^{-1/2},
\end{eqnarray*}
from inequality $|u_1\times \cdots  \times u_{j-1}-v_1 \times \cdots \times v_{j-1}|\le
|u_1-v_1|+\cdots+|u_{j-1}-v_{j-1}|$ , valid for complex numbers
$u_i,v_i$ with modulus less than $1$. Therefore, under the
different frames of dependence,
\begin{center}
\begin{tabular}{c|c|c}
Dependence & $|\cov(F_{t,n}(X_1,\ldots, X_{j-1}),G_{t,n}(X_{j})) |
$ & $T(k_n)$
 \\ \hline
$\theta$ & $\leq  \|t\|\cdot  k_n^{-1/2} \theta_{m_n}$& $\leq  \|t\|\cdot  k_n^{1/2} \theta_{m_n}$ \\
$\lambda$ & $\leq  \|t\|\cdot  k_n^{1/2} \lambda_{m_n}$& $\leq  \|t\|\cdot  k_n^{3/2} \lambda_{m_n}$ \\
\end{tabular}
\end{center}
and then, Theorem \ref{theolind} implies Proposition
\ref{subsampling}. $ \qquad\square$
\end{proof}
\begin{proof}
[Proof of Proposition \ref{kernel}]
Let $x\in \R$, define,
\begin{multline}
S_{n}=\sqrt{nh_n}\, \Big (\widehat f_X^{(n)}(x)-\E \widehat
f_X^{(n)}(x) \Big )=\sum_{i=1}^n Y_{i} \\
\mbox{where}~~~Y_i=\frac1{\sqrt{nh_n}}\Big (K\Big
(\frac{x-X_i}{h_n}\Big)-\E \Big (
K\Big(\frac{x-X_i}{h_n}\Big)\Big)\Big)=u(X_i)
\end{multline}
and the function $u$ depends also on $x$ and $n$. First, for
$\delta>0$,
\begin{eqnarray*}
A_n&=& \sum_{i=1}^n \E \big (\|Y_i\|^{2+\delta} \big ) \\
&= &(nh_n)^{-\delta/2} \cdot \frac 1 {h_n} \, \E \left ( \Big |
K\Big(\frac{x-X_i}{h_n}\Big)-\E \Big(
K\Big(\frac{x-X_i}{h_n}\Big)\Big) \Big |^{2+\delta} \right )\\
&\leq  &2 \cdot (nh_n)^{-\delta/2} \cdot \frac 1 {h_n} \,  \E
\left ( \Big |
K\Big(\frac{x-X_i}{h_n}\Big) \Big |^{2+\delta} \right )\\
&\leq  &2 \cdot (nh_n)^{-\delta/2} \cdot \|f_X\|_{\infty} \cdot
\int_\R  \big | K(s) \big |^{2+\delta}ds.
\end{eqnarray*}
(the boundedness of $K$ implies the convergence of the last
integral). As a consequence, $A_n \limiten 0$ when $nh_n\limiten
\infty$. Now,
\begin{eqnarray}\label{espe_kernel}
\nonumber  \left|\E (Y_i)\right|&\le& \|f_{X}\|_\infty\int_{\R}|u(s)|\,ds\\
\nonumber  &\le& \frac{2h_n}{\sqrt{nh_n}}\|f_X\|_\infty\int|K(v)|\, dv\\
&\le& C_1 \sqrt{\frac{h_n}{n}},\qquad \mbox{ for some constant }
C_1>0.
\end{eqnarray}
Moreover, using changes in variables $v=(x-s)/h_n$, $v'=(x-s')/h_n$,
\begin{eqnarray}\nonumber
\cov (Y_j,Y_i)&=&\int_{\R^2}u(s)u(s')(f_{j,i}(s,s')-f_X(s)f_X(s'))\, ds ds'\\
\nonumber
\left|\cov (Y_j,Y_i)\right|&\le& (\|f_{j,i}\|_\infty+\|f_X\|_\infty^2)\int_{\R^2}|u(s)||u(s')|\, dsds'\\
\nonumber
&\le& \frac{4h_n^2}{nh_n}\,(\|f_{j,i}\|_\infty+\|f_X\|_\infty^2)\left(\int|K(v)|\, dv\right)^2\\
\label{eqnoyaucov} &\le& C_2 \frac{h_n}{n},\qquad \mbox{ for some
constant } C_2>0.
\end{eqnarray}
The function $K$ is supposed to be a bounded and Lipschitz
function, the same for $u$,
$$
\|u\|_\infty\le 2\, \|K\|_\infty \cdot  \frac 1
{\sqrt{nh_n}}~~\mbox{and}~~\Lip u=2\,\lip K\cdot \frac 1
{h_n\sqrt{nh_n}}.
$$
Therefore, using the $\eta$-weak dependence inequality of the time
series $(X_i)_{i\in \N}$, with always $Y_i=u(X_i)$, there exists
$C>0$ such that,
\begin{eqnarray*}
|\cov (Y_0,Y_r)| \le C \cdot u_{n,r}~\mbox{with}~u_{n,r}=\left \{
\begin{array}{ll}
\displaystyle{\frac { \theta_r}{nh_n^2}}& \mbox{for the
$\theta$-weak dependence} \\
\displaystyle{\frac { \lambda_r}{nh_n^3}}& \mbox{for the
$\lambda$-weak dependence}
\end{array} \right .
\end{eqnarray*}
As a consequence of both (\ref{eqnoyaucov}) and the previous
inequalities, there exists constants $C_3>0$ such that
\begin{eqnarray}\label{cov_kernel}
|\cov (Y_0,Y_r)| \le   C_3 \cdot  \Big  ( \frac {h_n} n \wedge u_{n,r}
\Big ),\qquad \mbox{ for all } r\in \N.
\end{eqnarray}
Finally, we also quote that for $i \in \N$,
\begin{eqnarray*}
\v( Y_i)&= & \int_\R f_X(t) \cdot u^2(t)\,dt\\
& =& \frac 1 n \, \int_\R f_X(x-h_ns)\cdot K^2(s)\, ds.
\end{eqnarray*}
From the assumptions on functions $f_X$ and $K$, the Lebesgue
dominated convergence Theorem can be applied and therefore,
\begin{eqnarray}\label{var_kernel}
n \cdot \v( Y_i) \limiten  \int_\R  f_X(x)\cdot K^2(s)\, ds.
\end{eqnarray}
Using the relations (\ref{espe_kernel}),  (\ref{cov_kernel}) and
(\ref{var_kernel}), then,
\begin{eqnarray*}
\v (S_n)&=&n \cdot \v Y_0+2\sum_{i=1}^{n-1}(n-i)\cov(Y_0,Y_i)\\
\left|nh_n \cdot \v  \Big (\widehat f_X^{(n)}(x)\Big )-f_X(x)\int_\R
K^2(t)\,dt\right|&\le& o(1)+2C_3 \cdot \sum_{i=1}^{n-1}\Big  ( h_n
\wedge (n \cdot u_{n,r})\Big ).
\end{eqnarray*}
Under the assumptions of the Proposition, if the right term of the
forthcoming inequality (\ref{inegalite}) converges to $0$, the
right term of the previous inequality converges to $0$ and thus,
\begin{eqnarray}\label{var_f}
nh_n \cdot \v  \Big (\widehat f_X^{(n)}(x)\Big ) \limiten f(x)\int_\R
K^2(t)dt.
\end{eqnarray}
Now, we are going to bound $T(n)$ for applying Lemma \ref{lemlind'}.
Let $x \in \R$ and $t\in \R$. First we can write,
$$
T(n)=\sum_{j=1}^{n-1} \left |\cov(F_{x,t}(X_1,\ldots,
X_{j-1}),G_{x,t}(X_{j})) \right |,
$$
where $F_{x,t}(X_1,\ldots, X_{j-1})=\exp\Big (it(Y_1+\cdots+Y_{j-1})
\Big )$ and $G_{x,t}(X_{j})=\exp \big ( i t Y_j \big )$, with always
$Y_k=u(X_k)$. In order to compute a bound for $T(n)$ we need the
following decomposition due to Rio (2000),
$$
F_{x,t}(X_1,\ldots, X_{j-1})=\sum_{k=1}^{j-1}(e^{i t S_k}-e^{i t
S_{k-1}}),\quad \mbox{with}~~S_k=Y_1+\cdots+Y_k~~\mbox{and}~~S_0=0.
$$
Thus,
$$
\cov(F_{x,t}(X_1,\ldots,
X_{j-1}),G_{x,t}(X_{j}))=\sum_{k=1}^{j-1}\cov \Big (e^{i t
S_k}-e^{it S_{k-1}},e^{it Y_j}\Big ).
$$
Consider a r.v. $Y_j^*$ independent from $(Y_1,\ldots, Y_{k-1})$,
with the same distribution than $Y_j$. Then,
\begin{eqnarray*}
\Big |\cov \Big (e^{it S_k}-e^{i t S_{k-1}},e^{it Y_j}\Big )\Big
|&=&\Big |\E \Big (\big ( e^{i t S_k}-e^{it S_{k-1}}\big ) \big (e^{it Y_j}-e^{i t Y_j^*}\big ) \Big )\Big |\\
&\le& |t|^2  \E \big (|Y_k|\cdot (|Y_j|+|Y_j^*|) \big )~~\mbox{from inequality }~|e^{ia}-e^{ib}|\leq |b-a|\\
&\le& C  \frac{h_n}{n},\qquad \mbox{ for some constant } C>0,
\end{eqnarray*}
from relations (\ref{espe_kernel}) and (\ref{cov_kernel}). From
another hand, one can write,
$$
\cov \Big (e^{i t S_k},e^{i t Y_j}\Big )=\cov \Big
(g_1(X_1,\ldots,X_k),g_2(X_j) \Big ),
$$
with $\|g_2 \|_\infty = 1$ and $\displaystyle{\|g_1 \|_\infty \leq
|t| \| S_k-S_{k-1} \|_\infty \leq |t| \cdot \|u\|_\infty \leq 2|t|
\|K\|_\infty \cdot \frac 1 {\sqrt{nh_n}}}$ and,
\begin{eqnarray*}
\frac {\left
|e^{it(u(x_1)+\cdots+u(x_k))}-e^{it(u(y_1)+\cdots+u(y_k))}\right
|}{|x_1-y_1|+\cdots+|x_k-y_k|} & \leq & |t| \cdot \frac {\left
|(u(x_1)+\cdots+u(x_k))-(u(y_1)+\cdots+u(y_k))\right
|}{|x_1-y_1|+\cdots+|x_k-y_k|}
 \\
& \leq & |t|\cdot  \lip u \cdot k.
\end{eqnarray*}
As a consequence, $\displaystyle{~~ \lip g_1 \leq 8 |t| \cdot\lip K\cdot
\frac k {h_n\sqrt{nh_n}}~~~\mbox{and}~~\lip g_2 \leq 4 |t|\cdot\lip K\cdot
\frac 1 {h_n\sqrt{nh_n}}}$. Using these results and the weak
dependence property of $X$, there exists $C>0$ such that for
$\ell=j-k\in[1,j]$,
\begin{eqnarray*}
\Big |\cov \Big (e^{i t S_k}-e^{i t S_{k-1}},e^{i t Y_j}\Big )\Big
| &\le&C \cdot  \frac{h_n}{n}\wedge v_{n,k,\ell} \\
\mbox{with}~~~v_{n,k,\ell}&=&\left \{ \begin{array}{ll}
\displaystyle{k^2  \Big ( \frac{1}{nh_n^{3}}\vee
\frac{1}{n^{1/2}h_n^{3/2}}\Big )\cdot \lambda_{\ell} }& \mbox{for
the $\lambda$-weak dependence} \\
\displaystyle{\frac{1}{nh^2_n} \cdot \theta_{\ell} }& \mbox{for
the $\theta$-weak dependence}
\end{array}
\right . .
\end{eqnarray*}
This implies,
\begin{eqnarray}
\nonumber T(n)&\le&C  \sum_{j=1}^n\sum_{\ell=1}^j \frac{h_n}{n}
\wedge
v_{n,n,\ell} \\
\label{inegalite}
&\le & C  \sum_{\ell=1}^n {h_n}\wedge \big ( n \cdot v_{n,n,\ell}\big ) \\
\nonumber &\le & C \cdot h_n^{1-\beta} n^\beta  \sum_{\ell=1}^n
v_{n,n,\ell}^\beta
\end{eqnarray}
with $\beta \geq 0$ (analogously to the case of $\v\widehat
f^{(n)}_X(x)$). Since $nh_n\limiten \infty$ and $h_n \limiten 0$,
one obtains that $T(n)\limiten 0$ under the different conditions
satisfied by $h_n$ and the weak dependence sequence. Then, all the
conditions of Theorem \ref{theolind} are satisfied, which implies
Proposition \ref{kernel}. $ \qquad\square$
\end{proof}
\begin{proof}[Proof of Corollaries \ref{regu1} and \ref{regu2}]
Under the assumptions on $K$, from Prakasa Rao (1983),
$$
\E \Big (\widehat f_X^{(n)}(x)\Big )= \left \{
\begin{array}{ll}
\displaystyle{f_X(x)+h_n^p \cdot (1 +o(1))\cdot f_X^{(p)}(x)\frac
1 {p !} \, \int _\R t^{p} K(t)\, dt }&\mbox{if the regularity of
$f_X$ is $p \in
\N^*$} \\
\displaystyle{f_X(x)+O(h_n^{[\rho]})}&\mbox{if the regularity of
$f_X$ is $\rho \notin \N^*$}
\end{array}
\right . .
$$
It implies the optimal choice of convergence rate of $h_n$,
following this two cases, and the conditions on the convergence
rates of the different frames of weak dependent. $ \qquad\square$
\end{proof}
\begin{proof}[Proof of Proposition \ref{subkernel}]
This proof is
quite the same than the proof of Proposition \ref{kernel} and
therefore we omit the details. Let $x\in \R$, define,
\begin{multline}
S_{n}=\sqrt{k_nh_n}\, \Big (\widehat f_X^{(n,m_n)}(x)-\E \widehat
f_X^{(n,m_n)}(x) \Big )=\sum_{i=1}^{k_n} Y_{i} \\
\mbox{where}~~~Y_i=\frac1{\sqrt{k_nh_n}}\Big (K\Big
(\frac{x-X_{im_n}}{h_n}\Big)-\E \Big (
K\Big(\frac{x-X_{im_n}}{h_n}\Big)\Big)\Big)=u(X_{im_n}),
\end{multline}
and the function $u$ depends both on $x$ and $n$. First, for
$\delta>0$,
\begin{eqnarray*}
A^{(m_n)}_n&=& \sum_{i=1}^{k_n} \E \big (\|Y_i\|^{2+\delta} \big )
\leq  2 \cdot (k_nh_n)^{-\delta/2} \cdot \|f_X\|_{\infty} \cdot \int_\R \big |
K(s) \big |^{2+\delta}ds.
\end{eqnarray*}
As a consequence, $A^{(m_n)}_n \limiten 0$ when $k_nh_n\limiten
\infty$. Moreover, we have,
\begin{eqnarray*}\label{mom_subkernel}
\left|\E (Y_i)\right| &\le& C_1 \sqrt{\frac{h_n}{k_n}},\qquad
\mbox{for some constant } C_1>0;\\
|\cov (Y_0,Y_r)| &\le&   C_3 \cdot  \Big  ( \frac {h_n} {k_n} \wedge
u_{k_n,m_n \cdot r} \Big ),\qquad \mbox{for some constant } C_3>0
\mbox{ and all } r\in \N,
\end{eqnarray*}
with the sequence $(u_{p,q})$ defined in the proof of the
Proposition \ref{kernel}. We have also,
\begin{eqnarray*}
k_n \cdot \v( Y_i) \limiten  \int_\R  f_X(x)\cdot K^2(s)\, ds,
\end{eqnarray*}
and thus $\displaystyle{~~\left|k_nh_n \cdot \v  \Big (\widehat
f_X^{(n,m_n)}(x)\Big )-f_X(x)\int_\R K^2(t)\,dt\right| \le
o(1)+2C_3 \cdot \sum_{i=1}^{k_n-1}\Big  ( h_n \wedge (k_n \cdot u_{k_n,m_n
\cdot r})\Big )}$. Under the conditions on $h_n$ and $m_n$,
\begin{eqnarray}\label{var_subf}
k_nh_n \cdot \v  \Big (\widehat f_X^{(n,m_n)}(x)\Big ) \limiten
f(x)\int_\R K^2(t)dt.
\end{eqnarray}
For bounding $T(k_n)$, one writes again
$\mbox{with}~~S_k=Y_1+\cdots+Y_k~~\mbox{and}~~S_0=0$,
$$
T(k_n)=\sum_{j=1}^{k_n-1} \left |\sum_{k=1}^{j-1}\cov \Big (e^{i t
S_k}-e^{it S_{k-1}},e^{it Y_j}\Big ) \right |,
$$
Thanks to the inequality,
\begin{eqnarray*}
\Big |\cov \Big (e^{it S_k}-e^{i t S_{k-1}},e^{it Y_j}\Big )\Big |
&\le& C  \frac{h_n}{k_n},\qquad \mbox{ for some constant } C>0,
\end{eqnarray*}
and with the sequence $(v_{n,k,\ell})$ of the previous proof, this
implies,
\begin{eqnarray*}
T(k_n) &\le & C \cdot \sum_{i=1}^{k_n}\Big ( h_n \wedge  (k_n\cdot  v_{k_n,k_n,m_n\ell})\Big ) \\
&\le & C \cdot h_n^{1-\beta} k_n^\beta  \sum_{\ell=1}^{k_n}
v_{k_n,k_n,m_n\ell}^\beta
\end{eqnarray*}
with $\beta \geq 0$. With $h_n=n^{-h}$ and $m_n=n^m$, where $h,m
\in (0,1)$, the condition $T(k_n) \limiten 0$, which implies
(\ref{var_subf}) and therefore the central limit theorem, can be
obtained for different choice of $h$ and $m$. After computations,
the following graphs provide the zonas (depending also of the
decay rate $\theta$ or $\lambda$ of weak dependence property)
where $h$ and $m$ can be chosen,
\[
\epsfxsize 7cm \epsfysize 7cm \epsfbox{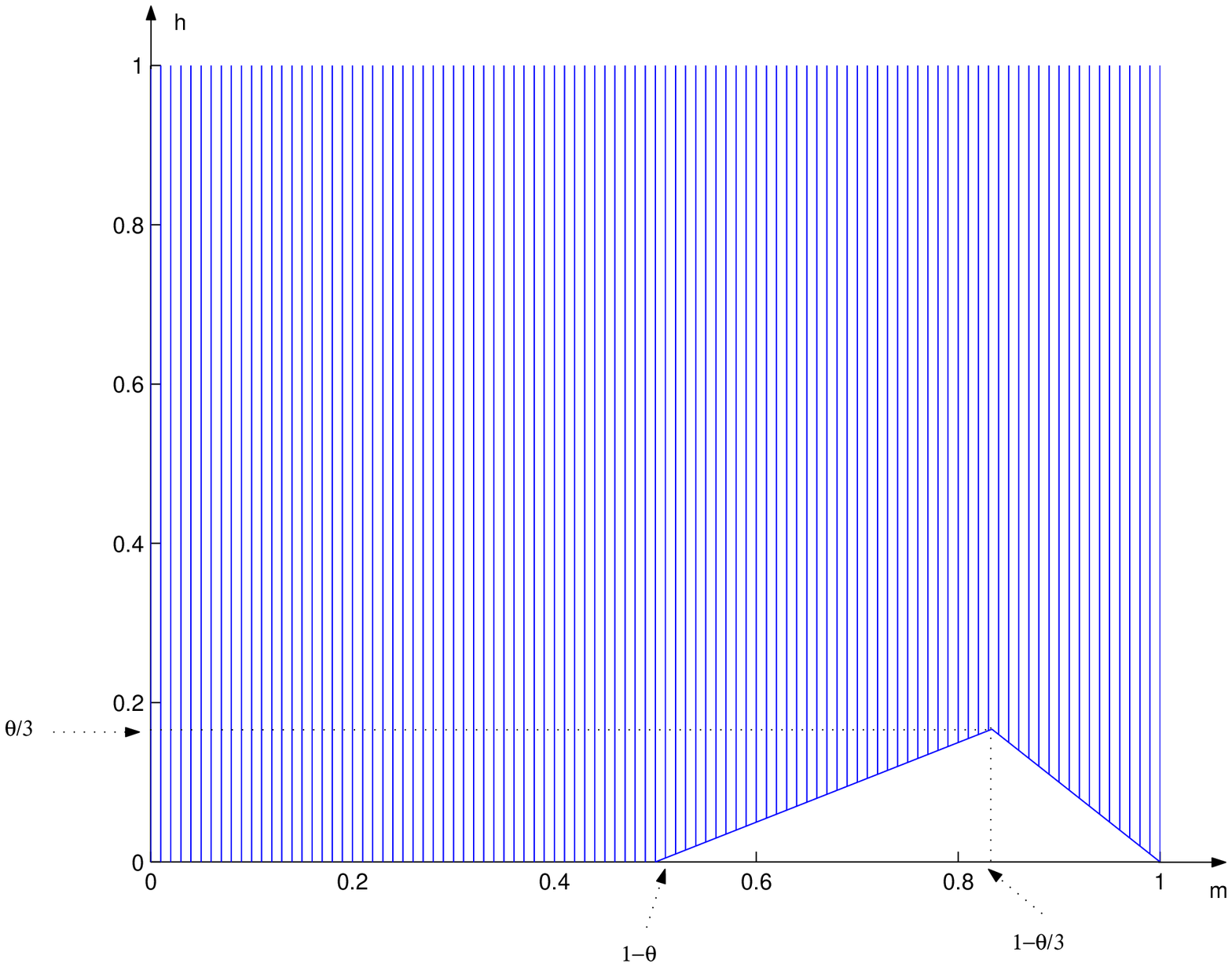}
\hspace*{0.7 cm} \epsfxsize 7cm \epsfysize 7cm
\epsfbox{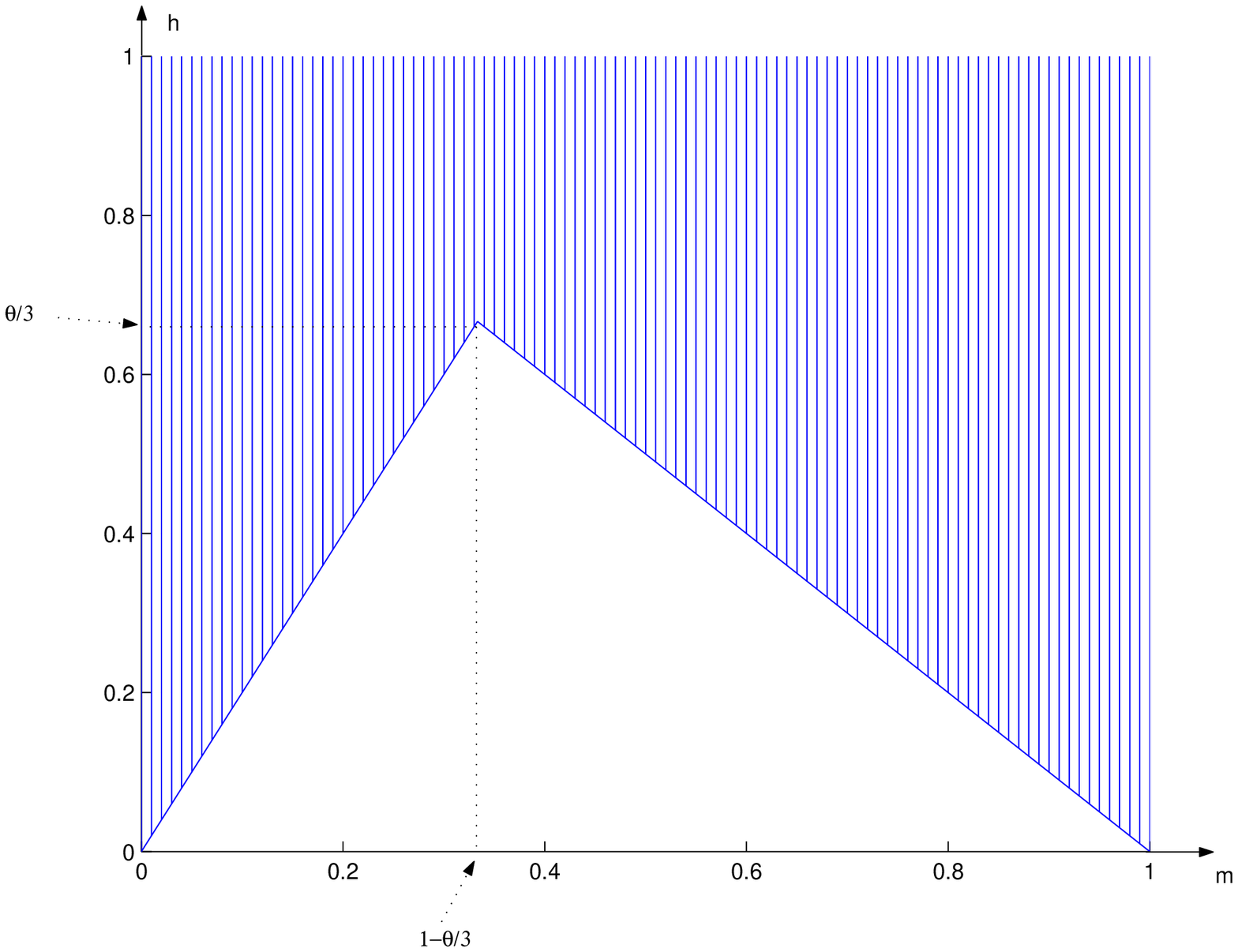}
\]
{\bf Figure 1 :} Conditions (the "white" zona) satisfied by
parameters $m$ and $h$ for obtaining the C.L.T. in the
$\theta$-weak dependence frame, when $0<\theta \leq 1$ (left) and
$1<\theta \leq 3$ (right).
\[
\epsfxsize 7cm \epsfysize 7cm \epsfbox{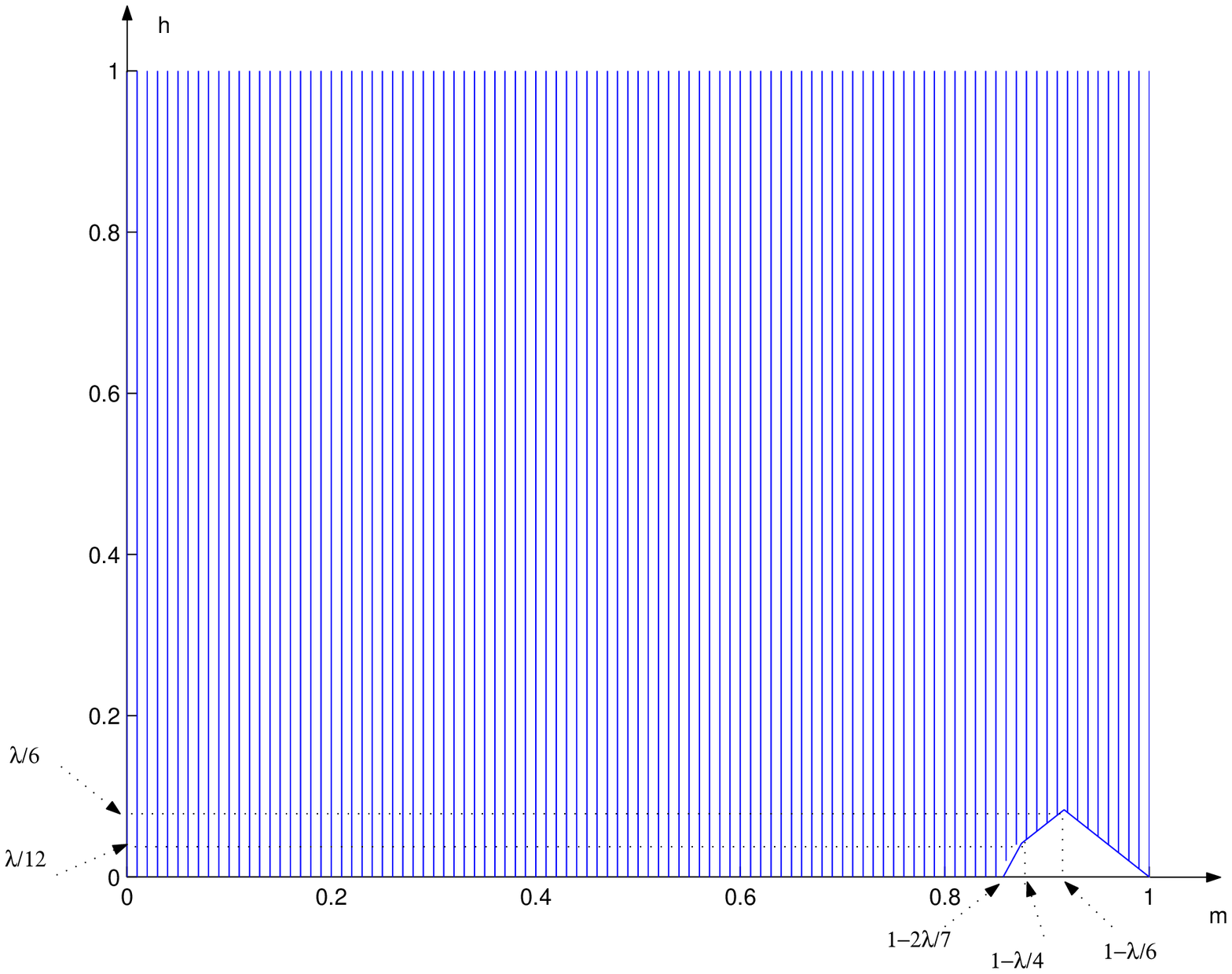}
\hspace*{0.7 cm} \epsfxsize 7cm \epsfysize 7cm
\epsfbox{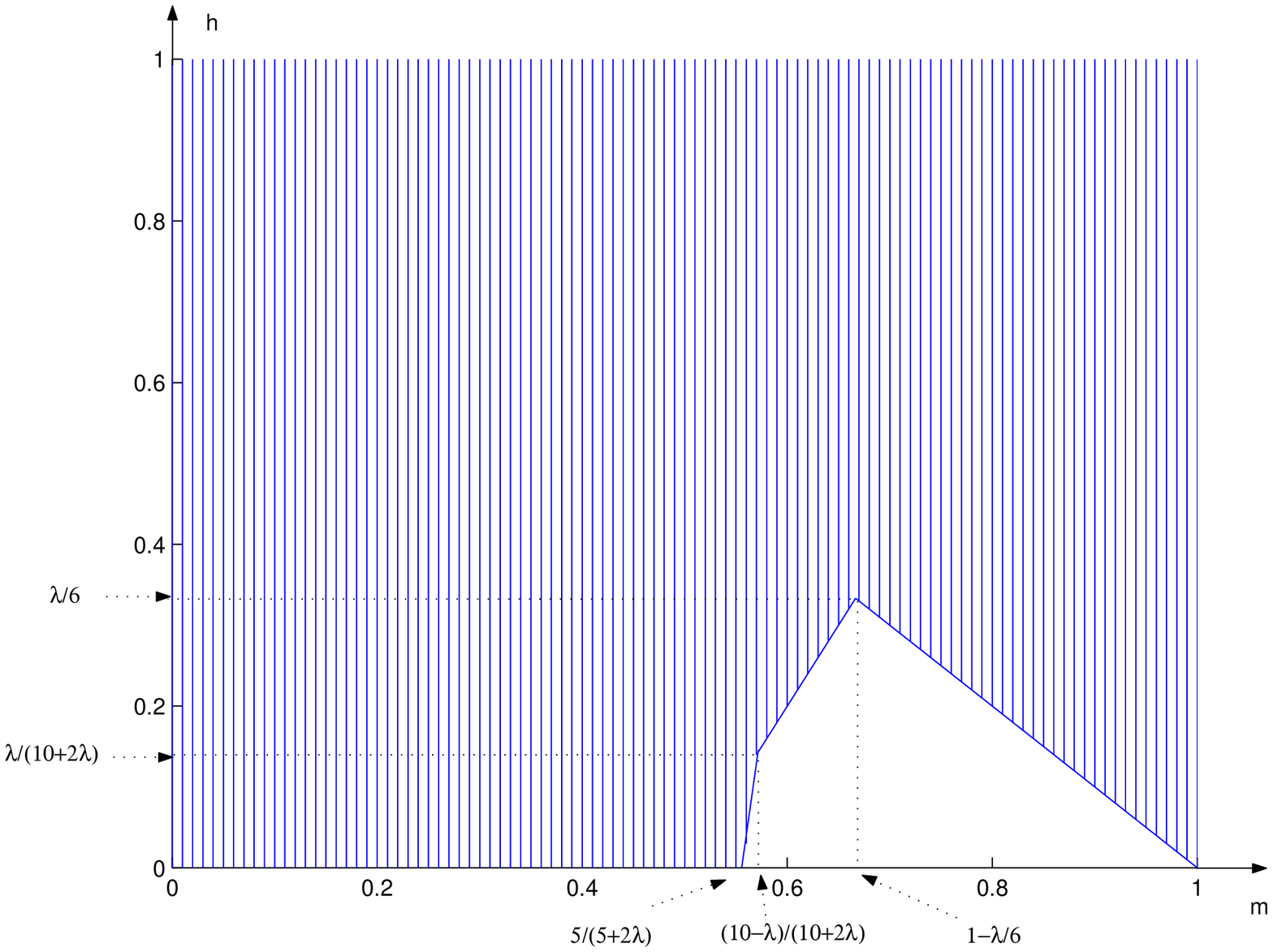}
\]
{\bf Figure 2 :} Conditions (the "white" zona) satisfied by
parameters $m$ and $h$ for obtaining the C.L.T. in the
$\lambda$-weak dependence frame, when $0<\lambda \leq 1$ (left)
and $1<\lambda\leq 6$ (right).
\\ ~\\
Finally, the "optimal" rate of convergence, in the sense of a
maximal $\sqrt{k_nh_n}=n^{(1-m-h)/2}$, is obtained from a
maximization of $1-m-h$. In every cases this occurs for the point
$(m,h)$ most below and left of the graph "white" zona. This
implies the optimal condition of the Proposition. $ \qquad\square$
\end{proof}

\end{document}